\documentclass[12pt]{article}

\usepackage{amsmath}
\usepackage{amsthm}
\usepackage{amsfonts}

\newenvironment{myproof}{\noindent {\it Proof} }{$\Box$ }
\newtheorem{theorem}{Theorem}[section]
\newtheorem{assertion}[theorem]{Proposition}
\newtheorem{lemma}[theorem]{Lemma}

\newtheorem{rem}[theorem]{Remark}

\newtheorem{example}[theorem]{Example}

\newcommand{\od}{\overset{d}{=}}

\newcommand{\zi}{Z_{\infty}}

\newcommand{\lin}{\underset{n\rightarrow\infty}{\lim}}

\newcommand{\mmp}{\mathbb{P}}
\newcommand{\me}{\mathbb{E}}

\newcommand{\mL}{\mathcal{L}}
\newcommand{\mr}{\mathbb{R}}

\newcommand{\omu}{\overline{\mu}}

\newcommand{\whQ}{\widehat{Q}}

\def\gl{\buildrel \rm def\over =}
\def\eqdist{\ {\buildrel d\over =}\ }
\def\1{{\bf 1}}

\begin{document}

\title{On distributional properties of perpetuities}\date{}
\author{Gerold Alsmeyer\footnote{e-mail
address: gerolda@math.uni-muenster.de}\\ \vspace{-.2cm}\noindent
\small{\emph{Institut f\"{u}r Mathematische Statistik},
\emph{Westf\"{a}lische Wilhelms-Universit\"{a}t M\"{u}nster}},\\
\small{\emph{Einsteinstra\ss e 62, 48149 M\"{u}nster, Germany}}\\
\\
Alex Iksanov\footnote {e-mail address:
iksan@unicyb.kiev.ua}\\ \vspace{-.2cm}\noindent
\small{\emph{Faculty of Cybernetics},
\emph{National T. Shevchenko University}},\\
\small{\emph{01033 Kiev, Ukraine}}\\ \\
Uwe R\"{o}sler\footnote{e-mail
address: roesler@math.uni-kiel.de}\\ \vspace{-.2cm}\noindent
\small{\emph{Mathematisches Seminar},
\emph{Christian-Albrechts-Universit\"at zu Kiel}},\\
\small{\emph{Ludewig-Meyn-Str.4, D-24098 Kiel, Germany}}}
\maketitle
\begin{quote}
\noindent {\footnotesize \textbf{SUMMARY.}} We study probability
distributions of convergent random series of a special structure,
called perpetuities. By giving a new argument,
we prove that such distributions are of pure type: degenerate,
absolutely continuous, or continuously singular. We further provide
necessary and sufficient criteria for the finiteness of $p$-moments, $p>0$
as well as exponential moments. In particular, a formula for the abscissa of convergence
of the moment generating function is provided. The results are illustrated with
a number of examples at the end of the article.
\medskip


\end{quote}

\section{Introduction and results}

Let $\{(M_n, Q_n): n=1,2,\ldots\}$ be a sequence of i.i.d.\ $\Bbb{R}^{2}$-valued
random vectors with generic copy $(M, Q)$. Put
$$\Pi_0\ \gl\ 1\quad\text{and}\quad\Pi_n\ \gl\ M_1 M_2 \cdots M_n,
\quad n=1,2,\ldots $$
Under conditions ensuring the almost sure convergence of the sequence
\begin{equation}\label{sequence}
Z_n\ \gl\ \sum_{k=1}^{n}\Pi_{k-1}Q_k,\quad n=1,2,\ldots
\end{equation}
the limiting random variable
\begin{equation}\label{per}
\zi\ \gl\ \sum_{k\ge 1}\Pi_{k-1}Q_k
\end{equation}
is often called perpetuity due to its occurrence in the realm of insurance and
finance as a sum of discounted payment streams. It also arises naturally from
shot-noise processes with exponentially decaying after-effect. In the study of
the (forward) iterated function system
\begin{equation}\label{ifsf}
\Phi_{n}\ \gl\ \Psi_{n}(\Phi_{n-1})\ =\ \Psi_{n}\circ...\circ\Psi_{1}
(\Phi_{0}),\quad n=1,2\ldots,
\end{equation}
where $\Psi_{n}(t)\gl Q_{n}+M_{n}t$ for $n=1,2,...$ and $\Phi_{0}$ is independent
of $\{(M_n, Q_n): n=1,2,\ldots\}$, the law of $\zi$ forms a
stationary distribution for this recursive Markov chain and thus a distributional
fixed point of the equation
\begin{equation}\label{fpeq}
\Phi\ \eqdist\ Q+M\Phi
\end{equation}
where as usual the variable $\Phi$ is assumed to be independent of $(M,Q)$.
Note that $\zi$ is obtained as the a.s.\ limit of the backward system associated
with (\ref{ifsf}) when started at $\Phi_{0}\equiv 0$, i.e.
\begin{equation}\label{ifsb}
\zi\ =\ \lim_{n\to\infty}\Psi_{1}\circ...\circ\Psi_{n}(0).
\end{equation}
This is in contrast to the forward sequence $\Psi_{n}\circ...\circ\Psi_{1}(0)$
which converges to $\zi$ in distribution only.

When focussing on equation (\ref{fpeq}), Vervaat \cite{Verv}
already showed that the law of $\zi$ forms the only possible solution, unless $Q+Mc=c$
a.s.\ for some $c\in\Bbb{R}$. Under the latter degeneracy condition, the
solutions to (\ref{fpeq}) are either all distributions on $\Bbb{R}$, or those symmetric
about $c$, or just the Dirac measure at $c$. This explains that subsequent work has
primarily dealt with a further study of the random variable $\zi$ and conditions that
ensure its existence. As to the last problem,
Goldie and Maller \cite{GolMal} gave the following complete characterization of
the a.s.\ convergence of the series in (\ref{per}). For $x>0$, define
\begin{equation}\label{A(x)}
A(x)\ \gl\ \int_0^x \mmp\{-\log |M|>y\}\ dy\ =\ \me\min(\log^{-}|M|,x).
\end{equation}
\begin{assertion}\label{exper} {\rm (\cite{GolMal}, Theorem 2.1)}
Suppose
\begin{equation}\label{nonzero}
\mmp\{M=0\}=0\quad\text{and}\quad\mmp\{Q=0\}<1.
\end{equation}
Then
\begin{equation}\label{cond2000}
\lin \Pi_n\,=\,0\ \text{a.s.}\quad\text{and}\quad I\ \gl\ \int_{(1,\infty)}
\dfrac{\log x}{A(\log x)}\ \mmp\{|Q|\in dx\}\,<\,\infty,
\end{equation}
and
\begin{equation}\label{conv}
\zi^{*}\ \gl\ \sum_{n\ge 1}|\Pi_{n-1}Q_n|\ <\ \infty\quad\text{a.s.}
\end{equation}
are equivalent conditions, which imply
\begin{equation*}
\lin Z_n\,=\,\zi\text{ a.s.\quad and}\quad|\zi|\ <\ \infty\quad\text{a.s.}
\end{equation*}
Moreover, if
\begin{equation}\label{degen}
\mmp\{Q+Mc=c\}<1\quad\text{for all } c\in \mr,
\end{equation}
and if at least one of the conditions in (\ref{cond2000}) fails to hold, then
$\lin |Z_n|=\infty$ in probability.
\end{assertion}

\begin{rem}\label{geom} \rm
As to condition (\ref{nonzero}) note that, if $\mmp\{Q=0\}=1$, $\zi$ trivially exists and
equals zero a.s. If $\mmp\{M=0\}>0$, then
$$ N\ \gl\ \min\{n\geq 1:M_n=0\}, $$
is a.s.\ finite and
$$ \zi\ =\ Z_n\ =\ \sum_{k=1}^N \Pi_{k-1}Q_k $$
for all $n\ge N$. Hence, in this case no condition on the distribution of $Q$
is needed to ensure the existence of $\zi$.
\end{rem}

Our first result states that the distribution of $\zi$ is pure.
\begin{theorem}\label{continuity}
If $\mmp\{M=0\}=0$ and $|\zi|<\infty$ a.s., then the distribution
of $\zi$ is either degenerate, absolutely continuous, or singular and
continuous.
\end{theorem}

To our knowledge this result is new in the given generality. It was obtained by
Grincevi\v{c}ius in \cite{Grinc74} under the additional condition $\me\log|M|\in (-\infty,
0)$. However, for the conclusion that $\zi$ is continuous if nondegenerate
his analytic argument is quite different from ours. This latter conclusion may also be
derived from his Theorem 1 in \cite{Grinc80}, as has been done in Lemma 2.1 in
\cite{Mal}. Let us further point out that $\zi=c$ a.s.\ for some $c\in\Bbb{R}$ implies
$\mmp\{Q+Mc=c\}=1$, as following from the fact that $\zi$ satisfies
(\ref{fpeq}). Hence the law of $\zi$ is continuous whenever (\ref{nonzero}) and (\ref{degen})
are assumed.

If $M$ and $Q$ are independent, Pakes \cite{Pakes83} provided sufficient
conditions for the absolute continuity of the distribution of $\zi$. His proof relies
heavily upon studying the behavior of corresponding characteristic functions and
used moment assumptions as an indispensable ingredient. Without such assumptions it
is not clear how absolute continuity of the law of $\zi$ may be derived via an analytic
approach.

\begin{theorem}\label{mom1}
Assuming (\ref{nonzero}) and (\ref{degen}), the following assertions are
equivalent for any $p>0$:
\begin{eqnarray}
&&\me |M|^p<1\quad\text{and}\quad\me |Q|^p<\infty,\label{cond2001}\\
&&\me |\zi|^p<\infty,\label{cond3000}\\
&&\me \zi^{*p}<\infty,\label{cond3001}
\end{eqnarray}
where $\zi^{*}$ is defined in (\ref{conv}).

\end{theorem}

Theorem \ref{mom1} seems to be new in the stated generality but was given as
Proposition 10.1 in \cite{Kel} for the case that $M, Q\geq 0$ (in which (\ref{cond3000})
and (\ref{cond3001}) are clearly identical).
If $p>1$, Vervaat \cite{Verv} proved that (\ref{cond2001}) implies
(\ref{cond3000}).

\begin{rem}\label{extra} \rm
It is not difficult to see that further conditions equivalent to those in the previous
theorem are given by
\begin{eqnarray}
&&\me\sup_{n\ge 1}|\Pi_{n-1}Q_{n}|^{p}<\infty,\label{cond3002}\\
&&\me\sup_{n\ge 1}|Z_{n}|^{p}<\infty,\label{cond3003}\\
&&\me\Bigg(\sum_{n\ge 1}\Pi_{n-1}^{2}Q_{n}^{2}\Bigg)^{p/2}<\infty,\label{cond3004}
\end{eqnarray}
where $Z_{n}$ is defined in (\ref{sequence}). Some comments regarding the proof can
be found in Remark \ref{comments} after the proof of Theorem \ref{mom1}.
\end{rem}

Given any real-valued random variable $Z$, let us define
$$ r(Z)\ \gl\ \sup\{r>0:\me e^{r|Z|}<\infty\}, $$
called the abscissa of convergence of the moment generating
function of $|Z|$. Note that $\me e^{r(Z)|Z|}$ may be finite or
infinite.

Our next two results provide complete information on how $r(\zi)$
relates to $r(Q)$. For convenience we distinguish the cases where
$\mmp\{|M|=1\}=0$ and $\mmp\{|M|=1\}\in (0,1)$. Recall that if
conditions (\ref{nonzero}) and (\ref{degen}) hold then the law of
$\zi$ is nondegenerate if $|\zi|<\infty$ a.s.

\begin{theorem}\label{perexp}
Suppose (\ref{nonzero}), (\ref{degen}) and $\mmp\{|M|=1\}=0$,
and let $s>0$. Then $\me e^{s|\zi|}<\infty$ holds if, and
only if,
\begin{eqnarray}
\mmp\{|M|<1\}=1\quad\text{and}\quad\me e^{s|Q|}<\infty.\label{cond4000}
\end{eqnarray}
In particular, if $\mmp\{|M|<1\}=1$, then $r(\zi)=r(Q)$.
\end{theorem}

\begin{theorem}\label{perexp2}
Suppose (\ref{nonzero}), (\ref{degen}) and $\mmp\{|M|=1\}\in (0,1)$, and let $s>0$.
Then $\me e^{s|\zi|}<\infty$ holds if,  and only if,
\begin{eqnarray}
&&\mmp\{|M|\leq 1\}=1,\quad \me e^{s|Q|}<\infty\quad\text{and}\label{cond4001}\\
&&b_{-}b_{+}<(1-a_{-})(1-a_{+}),\label{cond4002}
\end{eqnarray}
where $a_{\pm}=a_{\pm}(s)\gl\me e^{\pm sQ}\1_{\{M=1\}}$ and
$b_{\pm}=b_{\pm}(s)\gl \me e^{\pm sQ}\1_{\{M=-1\}}$. In
particular, if $\mmp\{|M|\leq 1\}=1$ and $\mmp\{|M|=1\}\in (0,1)$,
then $r(\zi)=\min\big(r(Q),r^{*}(M,Q)\big)$, where
$$ r^{*}(M,Q)\ \gl\ \sup\{r>0: b_{-}(r)b_{+}(r)<(1-a_{-}(r))(1-a_{+}(r))\}. $$
\end{theorem}

The reader should notice $\max(a_{-}+b_{-},a_{+}+b_{+})<1$ provides a sufficient
condition for (\ref{cond4002}).

We are aware of two papers that deal with the existence of exponential
moments of $|\zi|$. By using the contraction principle, Goldie and Gr\"{u}bel
\cite{GolGru} proved that, if $|M|\leq 1$ a.s. and $\me e^{r|Q|}<\infty$ for
some $r>0$, then $\me e^{t |\zi|}<\infty$ for $0\leq t< \sup\{\theta: \me
e^{\theta|Q|}|M|<1\}$ (see their Theorem 2.1). For the case of nonnegative
$M$ and $Q$, a stronger result was obtained by Kellerer \cite{Kel}, namely
that $\me e^{r\zi}<\infty$ for some $r>0$ iff $M\leq 1$ a.s. and $\me
e^{tQ}<\infty$ for some $t>0$ (see his Proposition 10.2).

The rest of the paper is organized as follows. All proofs are
given in Section 2. Section 3 collects a number of examples which
illustrate our main results, followed by a number of concluding remarks
in Section 4.

\section{Proofs}

Let us start by pointing out that, if $|\zi|<\infty$ a.s.,
\begin{equation}\label{pereqdistr}
\zi\ =\ Q_{1}+M_{1}\zi^{(1)}\ =\ Q^{(m)}+\Pi_{m}\zi^{(m)},
\end{equation}
holds true for each $m\ge 1$, where (setting $\Pi_{k:l}\gl M_{k}\cdot...\cdot M_{l}$)
\begin{equation}\label{iterate}
Q^{(m)}\ \gl\ \sum_{k=1}^{m}\Pi_{k-1}Q_{k}\quad\text{and}\quad
\zi^{(m)}\ \gl\ Q_{m+1}+\sum_{k\ge m+2}\Pi_{m+1:k-1}Q_k.
\end{equation}
The latter variable is a copy of $\zi$ and independent of
$(M_1,Q_1),...,(M_{m},Q_{m})$. We thus see that $\zi$ may be viewed as the perpetuity
generated by i.i.d.\ copies of $(\Pi_{m},Q^{(m)})$ for any fixed $m\ge 1$.
\bigskip

\begin{myproof}{\emph{of Theorem \ref{continuity}}}.
Suppose immediately that $\mmp\{Q=0\}<1$ and (\ref{degen}) hold true, for otherwise the
law of $\zi$ is clearly degenerate. By Proposition \ref{exper} we thus also have
that $\Pi_{n}\to 0$ a.s.

We first show that if the law of $\zi$ is having atoms then again it must be degenerate.
Let $b_{1},...,b_{d}$ denote the atoms with maximal probability $\varrho$,
say. Notice that $d\leq\varrho^{-1}$. In view of (\ref{pereqdistr}) we have
\begin{equation}\label{cont1}
\mmp\{\zi=b_i\}\ =\ \sum_{a\in A}\mmp\{Q^{(m)}+\Pi_{m}a=b_i\}\,\mmp\{\zi=a\},\quad
i=1,...,d
\end{equation}
for each $m=1,2,...$,
where $A$ is the set of all atoms of the distribution of $\zi$.
Since $\mmp\{M=0\}=0$, we have $\sum_{a\in A}\mmp\{Q^{(m)}+\Pi_{m}a=b_i\}\leq 1$.
Now use $\mmp\{\zi=a\}\leq\mmp\{\zi=b_i\}$ to conclude that
(\ref{cont1}) can only hold if the summation extends only over $b_j$,
$j=1,...,d$, and so
\begin{equation}\label{cont2}
\sum_{j=1}^{d}\mmp\{Q^{(m)}+\Pi_{m}b_j=b_i\}\ =\ 1,\quad i=1,...,d,
\end{equation}
for each $m=1,2,...$ By letting $m$ tend to infinity and using $(\Pi_{m},Q^{(m)})
\to (0,\zi)$ a.s. in (\ref{cont2}), we arrive at
\begin{equation}\label{laplace}
\mmp\{\zi=b_i\}\ =\ d^{-1},\quad i=1,...,d.
\end{equation}
In order to see that this already yields degeneracy of $\zi$, suppose $d\ge 2$ and let
$U,V$ be independent copies of $\zi$ which are also independent of
$\{(M_{n},Q_{n}):n=1,2,...\}$. Put
$\zi^{(s)}\gl U-V$, clearly a symmetrization of $\zi$ with support given as
$\Gamma\gl\{b_{i}-b_{j}: i,j=1,...,d\}$. Since $Q^{(m)}+\Pi_{m}U\od
Q^{(m)}+\Pi_{m}V\od\zi$ for each $m=1,2,...$, we see
that
$$ D_{m}\ \gl\ \big(Q^{(m)}+\Pi_{m}U\big)-\big(Q^{(m)}+\Pi_{m}V\big)\ =\ \Pi_{m}\zi
^{(s)} $$
has a support $\Gamma_{m}$ contained in $\Gamma$. Put $\gamma_{*}\gl\min(\Gamma\cap
(0,\infty))$ and $\gamma^{*}\gl\max\Gamma$.
Using the independence of
$\Pi_{m}$ and $\zi^{(s)}$ in combination with $\mmp\{M=0\}=0$, we now infer
\begin{eqnarray*}
0 &=& \mmp\{|D_{m}|\in (0,\gamma_{*})\}\ =\ \mmp\{|\Pi_{m}\zi^{(s)}|\in (0,\gamma_
{*})\}\\
&\ge& \mmp\{|\Pi_{m}|<\gamma_{*}/\gamma^{*})\,\mmp\{|\zi^{(s)}|\in (0,\gamma^{*}]\}
\end{eqnarray*}
and therefore $\mmp\{|\Pi_{m}|<\gamma_{*}/\gamma^{*})=0$ because
$\mmp\{|\zi^{(s)}|\in (0,\gamma^{*}]\}=1-d\varrho^{2}>0$. But this contradicts
$\Pi_{m}\to 0$ a.s.\ and so $d=1$, i.e.\ $\zi=b_{1}$ a.s. by (\ref{laplace}).

It remains to verify that a continuously distributed $\zi$ is of pure type.
Apart from minor modifications, the following argument is due to
Grincevi$\check{c}$ius \cite{Grinc74} and stated here for completeness.

Let $\phi(t)$ be the characteristic function (ch.f.) of $\zi$. By Lebesgue's\nobreak\
decomposition theorem $\phi=\alpha_{1}\phi_{1}+\alpha_{2}\phi_{2}$, where
$\alpha_{1},\alpha_{2}\ge 0$, $\alpha_{1}+\alpha_{2}=1$, and $\phi_{1},\phi_{2}$ are
the ch.f.\ of the absolutely continuous and the continuously singular components of the
law of $\zi$, respectively. Suppose $\alpha_{1}>0$ so that $\phi=\phi_{1}$ must be
verified.

Since the law of $\zi$ satisfies the stochastic fixed point
equation (\ref{fpeq}), we infer in terms of its ch.f.
\begin{equation}\label{chffp}
\phi(t)\ =\ \me e^{itQ}\phi(Mt)
\end{equation}
and thus
\begin{equation}\label{Lebesgue}
\alpha_{1}\phi_{1}(t)+\alpha_{2}\phi_{2}(t)\ =\ \alpha_{1}\me e^{itQ}\phi_{1}
(Mt)+\alpha_{2}\me e^{itQ}\phi_{2}(Mt).
\end{equation}
It is easily seen that $\me e^{itQ}\phi_{1}(Mt)$ is the ch.f.\ of an absolutely
continuous distribution, because this is true for $\phi_{1}$ and $\mmp\{M=0\}=0$. If
$$ \alpha_{2}\me e^{itQ}\phi_{2}(Mt)=\alpha_{3}\phi_{3}(t)+\alpha_{4}\phi_{4}(t), $$
where $\alpha_{3},\alpha_{4}\ge 0$, $\alpha_{3}+\alpha_{4}=1$, and
$\phi_{3},\phi_{4}$ are the absolutely continuous and the continuously singular
components, respectively, then the uniqueness of the Lebesgue decomposition renders in
(\ref{Lebesgue}) that $\alpha_{1}\phi_{1}(t)=\alpha_{1}\me e^{itQ}\phi_{1}
(Mt)+\alpha_{2}\alpha_{3}\phi_{3}(t)$ and thus upon setting $t=0$ that
$\alpha_{2}\alpha_{3}=0$. Consequently, $\phi_{1}(t)=\me e^{itQ}\phi_{1}(Mt)$
which means that $\phi_{1}$ is also a solution to the functional equation (\ref{chffp}).
By considering the bounded continuous function $\phi-\phi_{1}$ and utilizing
$\phi(0)-\phi_{1}(0)=0$ in combination with $\Pi_{n}\to 0$ a.s., we infer upon iterating
(\ref{chffp}) for $\phi-\phi_{1}$ and an appeal to the dominated convergence theorem
that
$$ |\phi(t)-\phi_{1}(t)|\ \le\ \lim_{n\to\infty}\me\big|\phi(\Pi_{n}t)-
\phi(\Pi_{n}t)\big|\ =\ 0 $$
for all $t\ne 0$ and so $\phi=\phi_{1}$.
\end{myproof}

\bigskip
\begin{myproof}\emph{of Theorem \ref{mom1}}.
Since $|\zi|\le\zi^{*}$ it suffices to prove
"(\ref{cond2001})$\Rightarrow$(\ref{cond3001})" and
"(\ref{cond3000})$\Rightarrow$(\ref{cond2001})". A proof of the first implication,
though known from \cite{Kel}, is easy and thus stated for completeness here.
\medskip

"(\ref{cond2001})$\Rightarrow$(\ref{cond3001})":
If $0<p\le 1$, just use the subadditivity of $x\mapsto x^{p}$ on $[0,\infty)$ in
combination with the independence of $\Pi_{k-1}$ and $Q_{k}$ for each $k$ to infer
\begin{eqnarray*}
\me\zi^{*p} &\le& \sum_{k\ge 1}\me|\Pi_{k-1}|^{p}\,\me|Q_{k}|^{p}\\
&=& \sum_{k\ge 1}(\me|M|^{p})^{k-1}\me|Q|^{p}
\ =\ \dfrac{\me|Q|^{p}}{1-\me|M|^{p}}\ <\ \infty.\\
\end{eqnarray*}
If $p>1$, a similar inequality holds for $\|\zi\|_{p}$, where $\|\cdot\|_{p}$ denotes the
usual $L_{p}$-norm. Namely, by Minkowski's inequality,
\begin{eqnarray*}
\|\zi^{*}\|_{p}\ \le\ \sum_{k\ge 1}\|\Pi_{k-1}Q_{k}\|_{p}\ =\ \sum_{k\ge 1}\|M\|_{p}
^{k-1}\|Q\|_{p}\ =\ \dfrac{\|Q\|_{p}}{1-\|M\|_{p}}\ <\ \infty.
\end{eqnarray*}
\medskip

"(\ref{cond3000})$\Rightarrow$(\ref{cond2001})" Let us start by pointing out that
(\ref{cond2001}) is equivalent to
\begin{equation}\label{cond2002}
\me|Q_{1}+M_{1}Q_{2}|^{p}<\infty\quad\text{and}\quad\me|M_{1}M_{2}|^{p}<1.
\end{equation}
which, in the notation introduced at the beginning of this section, is nothing but
condition (\ref{cond2001}) for the pair $(\Pi_{2},Q^{(2)})$. We only remark concerning
"(\ref{cond2002})$\Rightarrow$(\ref{cond2001})" that in the case $p\ge 1$, by Minkowski's
inequality,
\begin{eqnarray*}
\|Q_{1}\1_{\{|Q_{1}|\le b,|Q_{2}|\le c\}}\|_{p} &\le&
\|(Q_{1}+M_{1}Q_{2})\1_{\{|Q_{1}|\le b,|Q_{2}|\le c\}}\|_{p}\\
&+&\|M_{1}\1_{\{|Q_{1}|\le b\}}\|_{p}\,\|Q_{2}\1_{\{|Q_{2}|\le c\}}\|_{p}
\end{eqnarray*}
for all $b,c>0$ and therefore (upon letting $b$ tend to $\infty$ and picking $c$ large
enough)
$$ \|Q\|_{p}\ \le\ {\|Q_{1}+M_{1}Q_{2}\|_{p}+c\,\|M\|_{p}\over \mmp\{|Q|\le c\}^{1/p}}
\ <\ \infty. $$
If $0<p<1$ a similar argument using the subadditivity of $0\le x\mapsto x^{p}$ yields
the conclusion.
Next, we note that the conditional law of $Q_{1}+M_{1}Q_{2}$ given $\Pi_{2}$ cannot
be degenerate, for otherwise either $Q+cM=c$ or $(M_{1},Q_{1})=(1,c)$ a.s.\ for some
$c\in\Bbb{R}$ by Proposition 1 in \cite{Grinc81}. But both alternatives are here
impossible, the first by our assumption (\ref{degen}), the second by $|\zi|<\infty$ a.s.
Let us also mention that $|\zi|<\infty$ a.s.\ in combination with (\ref{degen})
ensures $\Pi_{n}\to 0$ a.s. by Theorem \ref{exper}.

The following argument based upon conditional symmetrization may be viewed as a
streamlined version of a similar one given in the proof of Proposition 3 in
\cite{IksRos}. Put $Q_{n}^{(2)}\gl Q_{2n-1}+M_{2n-1}Q_{2n}$ for $n=1,2,...$ and note that
$\{(M_{2n-1}M_{2n},Q_{n}^{(2)}):n=1,2,...\}$ is a family of independent copies of
$(\Pi_{2}, Q^{(2)})$. Let $\overline{Q}_{n}^{(2)}$ be a conditional symmetrization of
$Q_{n}^{(2)}$ given $M_{2n-1}M_{2n}$ such that
$(M_{2n-1}M_{2n},\overline{Q}_{n}^{(2)}),\,n=1,2,...$ are also i.i.d.
More precisely, $\overline{Q}_{n}^{(2)}=Q_{n}^{(2)}-\whQ_{n}^{(2)}$,
where $\{(M_{2n-1}M_{2n},Q_{n}^{(2)},\whQ_{n}^{(2)}):n=1,2,...\}$ consists of
i.i.d.\ random variables and $Q_{n}^{(2)},\whQ_{n}^{(2)}$ are conditionally
i.i.d.\ given $M_{2n-1}M_{2n}$.
By what has been pointed out above, the law of $\overline Q_{n}^{(2)}$, and thus also of
${Q}_{n}^{(2)}$, is\nobreak\ non\-degenerate.
Putting ${\cal B}_{n}\gl\sigma(M_{1},...,M_{n})$ for $n=1,2,...$, we now infer with the
help of L\'evy's symmetrization inequality (see \cite{Chow}, Corollary 5 on p.\ 72)
\begin{eqnarray*}
\mmp\Big(\max_{1\le k\le n}|\Pi_{2k-2}\overline{Q}_{k}^{(2)}|>x\Big|{\cal B}_{2n}\Big)
&\le& 2\,\mmp\Bigg(\bigg|\sum_{k=1}^{n}\Pi_{2k-2}\overline{Q}_{k}^{(2)}\bigg|>x
\bigg|{\cal B}_{2n}\Bigg)\\
&\le& 4\,\mmp\Bigg(\bigg|\sum_{k=1}^{n}\Pi_{2k-2}Q_{k}^{(2)}\bigg|>{x\over 2}\bigg|
{\cal B}_{2n}\Bigg)\\
&=& 4\,\mmp(|Z_{2n}|>x/2|{\cal B}_{2n})\quad\text{a.s.}
\end{eqnarray*}
for all $x>0$ and thus (recalling that the law of $\zi$ is continuous in the present situation
as pointed out right after Theorem \ref{continuity})
\begin{equation}\label{symm}
\mmp\Big\{\sup_{k\ge 1}|\Pi_{2k-2}\overline{Q}_{k}^{(2)}|>x\Big\}
\ \le\ 4\,\mmp\{|\zi|>x/2\}.
\end{equation}
As a consequence of this in combination with $\me|\zi|^{p}<\infty$ we conclude
$$ \me\sup_{k\ge 1}|\Pi_{2k-2}\overline{Q}
_{k}^{(2)}|^{p}\ \le\ 8\,\me|\zi|^{p}\ <\ \infty. $$
Now put $S_{0}\gl 0$ and
$$ S_{n} \gl \log|\Pi_{2n}|=\sum_{k=1}^{n}\log|M_{2k-1}M_{2k}|\quad\text{and}
\quad Y_{n} \gl \1_{\{\overline{Q}_{n}^{(2)}\ne 0\}}
\log|\overline{Q}_{n}^{(2)}| $$
for $n=1,2,...$ Then
$\{S_{n}:n=0,1,...\}$ forms an ordinary zero-delayed random walk with $S_{n}\to
-\infty$ a.s. (recall $\Pi_{n}\to 0$ a.s.\ from above), and
$$ \mmp\{Y_{n}=0\hbox{ i.o}\}\ =\ \mmp\big\{\overline{Q}_{n}^{(2)}\in\{0,1\}
\hbox{ i.o}\big\}\ =\ 0 $$
by the nondegeneracy and symmetry of the $\overline{Q}_{n}^{(2)}$. With this we see that
$$ \me\sup_{k\ge 1}|\Pi_{2k-2}\overline{Q}_{k}^{(2)}|^{p}\ =\ \me\exp\Big(
p\sup_{n\ge 0}(S_{n}+Y_{n+1})\Big)\ <\ \infty. $$
Since the pairs $(\log|M_{2n-1}M_{2n}|,Y_{n})$, $n=1,2,...$, are i.i.d., an
application of Lemma \ref{Glynn} stated below yields
$$ \me\sup_{k\ge 0}|\Pi_{2k}|^{p}\ =\ \me\exp\Big(p\sup_{k\ge 0}S_{k}\Big)
\ <\ \infty. $$
%(and then also $\me\sup_{k\ge 0}|\Pi_{k}|^{p}<\infty$ as one can easily check).
But we further have that $W\gl\sup_{k\ge 0}|\Pi_{2k}|^{p}$ and its copy $W'\gl
\sup_{k\ge 1}|\Pi_{3:2k}| ^{p}$ (setting $\Pi_{3:2}\gl 1$) satisfy
\begin{equation}\label{meanone}
W\ =\ \max(1,|\Pi_{2}|^{p}W')\ \ge\ |\Pi_{2}|^{p}W'\quad\text{and}\quad
\me W\ \ge\ \me|\Pi_{2}|^{p}\,\me W'
\end{equation}
whence $\me|\Pi_{2}|^{p}=\me|M_{1}M_{2}|^{p}\le 1$. In order to conclude strict
inequality note first that $\me|\Pi_{2}|^{p}=1$ in (\ref{meanone}) would give
$|\Pi_{2}|^{p}W'=W\ge 1$ a.s. But since $p\gl\mmp\{W=1\}\ge\mmp\{\sup_{n\ge 1}|\Pi_{2n}|
<1\}>0$ as argued below, the independence of
$W'$ and $\Pi_{2}$ would further imply
\begin{eqnarray*}
\mmp\{|\Pi_{2}|^{p}W'<1\}\ \ge\ \mmp\{|\Pi_{2}|<1,W'=1\}\ =\ p\,\mmp\{|\Pi_{2}|<1\}
\ >\ 0
\end{eqnarray*}
which is a contradiction. Therefore
$\me|\Pi_{2}|^{p}<1$, that is the second half of (\ref{cond2002}) holds true.

In order to show $\mmp\{\sup_{n\ge 1}|\Pi_{2n}|<1\}>0$, let us recall that
$S_{n}=-\log|\Pi_{2n}|$, $n=0,1,...$, forms an ordinary random walk converging a.s.\ to
$-\infty$ (as $\Pi_{n}\to 0$ a.s.). Consequently, the associated first strictly
descending ladder epoch $\tau_{-}\gl\inf\{n:S_{n}<0\}$ has finite mean (see Cor.\ 1
on p.\ 153 in \cite{Chow}) and with its dual $\tau_{+}\gl\inf\{n:S_{n}\ge 0\}$ it
satisfies the relation
$$ \mmp\{\tau_{+}=\infty\}\ =\ 1/\me\tau_{-}\ >\ 0, $$
see Thm.\ 2 on p.\ 151 in\cite {Chow}. But $\{\tau_{+}=\infty\}=\{\sup_{n\ge
1}|\Pi_{2n}|<1\}$.

Left with the first half of (\ref{cond2002}), namely
$\|Q^{(2)}\|_{p}<\infty$, use (\ref{pereqdistr}) with $m=2$ rendering
$|Q^{(2)}|\le |\zi|+|\Pi_{2}\zi^{(2)}|$ and therefore
$$ \|Q^{(2)}\|_{p}\ \le\ \|\zi\|_{p}(1+\|\Pi_{2}\|_{p})\ <\ \infty $$
in the case $p\ge 1$. The case $0<p<1$ is handled similarly.
This completes the proof of the theorem.
\end{myproof}

\begin{rem}\label{comments} \rm Here are a few comments
on how to obtain the additional equivalences stated in Remark \ref{extra}:\newline
"(\ref{cond3002})$\Rightarrow$(\ref{cond2001})" is contained in the above proof of
Theorem
\ref{mom1}.\newline
"(\ref{cond2001})$\Rightarrow$(\ref{cond3002})": Use
"(\ref{cond2001})$\Rightarrow$(\ref{cond3001})" and $\sup_{n\ge 1}|\Pi_{n-1}Q_{n}|
\le\zi^{*}$.\newline
"(\ref{cond2001})$\Rightarrow$(\ref{cond3003})": Use
"(\ref{cond2001})$\Rightarrow$(\ref{cond3001})" and $\sup_{n\ge 1}|Z_{n}|\le\zi^{*}$.
\newline
"(\ref{cond3003})$\Rightarrow$(\ref{cond2001})": Use $|\zi|\le\sup_{n\ge 0}|Z_{n}|$
and then "(\ref{cond3000})$\Rightarrow$(\ref{cond2001})".\newline
"(\ref{cond3004})$\Rightarrow$(\ref{cond2001})": Use
$\sup_{n\ge 1}|\Pi_{n-1}Q_{n}|\le\big(\sum_{n\ge 1}\Pi_{n-1}^{2}Q_{n}^{2}
\big)^{1/2}$ and "(\ref{cond3002})$\Rightarrow$(\ref{cond2001})".\newline
"(\ref{cond2001})$\Rightarrow$(\ref{cond3004})": Use
"(\ref{cond2001})$\Rightarrow$(\ref{cond3001})" and $\big(\sum_{n\ge 1}\Pi_{n-1}^{2}
Q_{n}^{2}\big)^{1/2}\le\zi^{*}$.
\end{rem}

\bigskip
We continue with the lemma used at the end of the previous proof. It
contains a tail inequality first given in \cite{Goldie}. A similar result
was stated as Lemma 2 in \cite{IksRos}, but that result is
correct only for the case of independent $M$ and $Q$.

\begin{lemma}\label{Glynn}
Let $\{(X_k,Y_k): k=1,2,\ldots\}$ be a family of i.i.d.\ $\Bbb{R}^{2}$-valued random
vectors. Put $S_n\gl X_1+\cdots+X_n$ for $n=0,1,...$, $\xi\gl\sup_{n\ge 0}S_{n}$
and $\zeta\gl\sup_{n\ge 0}(S_{n}+Y_{n+1})$. Then
\begin{equation}\label{tailin}
\mmp\{\zeta>x\}\ \ge\ \mmp\{Y_{1}>y\}\,\mmp\{\xi>x-y\}
\end{equation}
for all $x,y\in\Bbb{R}$. Furthermore, if $\Phi:[0,\infty)\to [0,\infty)$ is any
nondecreasing, differentiable function, then
\begin{equation}\label{momentin}
\me\Phi(\xi)\ \le\ \Phi(0)+c\,\me\Phi(c+\zeta^{+})
\end{equation}
for a constant $c\in (1,\infty)$ that does not depend on $\Phi$.
\end{lemma}

\begin{myproof}
For any fixed $x,y\in\Bbb{R}$, put
$\tau\gl\inf\{k\geq 0: S_k> x-y\}$
with the usual convention $\inf\emptyset\gl\infty$. Note that $\{\xi>x-y\}=
\{\tau<\infty\}$ and $\{\zeta>x\}\supset\{\tau<\infty,Y_{\tau+1}>y\}$.
Inequality (\ref{tailin}) now follows from
\begin{eqnarray*}
\mmp\{\zeta>x\}&\ge&\mmp\{\tau<\infty,Y_{\tau+1}>y\}
\ =\ \sum_{n\ge 0}\mmp\{\tau=n,\,Y_{n+1}>y\}\\
&=& \mmp\{Y_{1}>y\}\sum_{n\ge 0}\mmp\{\tau=n\}\ =\ \mmp\{Y_{1}>y\}
\,\mmp\{\tau<\infty\}\\
&=& \mmp\{Y_{1}>y\}\,\mmp\{\xi>x-y\}.
\end{eqnarray*}
In order to get (\ref{momentin})
fix any $c>1$ such that $\mmp\{Y_{1}>-c\}\ge 1/c$.
Then (\ref{tailin}) with $y=-c$ provides us with
$$ \mmp\{\zeta+c>x\}\ \ge\ \mmp\{\xi>x\}/c $$
for $x\in\Bbb{R}$ which in combination with
$\me\Phi(\xi)-\Phi(0)=\int_{0}^{\infty}\Phi'(x)\,\mmp\{\xi>x\}\,dx$ finally shows
(\ref{momentin}).
\end{myproof}
\bigskip

\begin{myproof}{\emph{of Theorem \ref{perexp}}}. Let us define $\psi(t)\gl\me e^{t\zi}$,
$\widehat\psi(t)\gl\me e^{t|\zi|}$, $\varphi(t)\gl\me e^{tQ}$ and $\widehat{\varphi}
(t)\gl\me e^{t|Q|}$.
Note that $\psi(t)\le\widehat{\psi}(s)$ for all $t\in [-s,s]$ and $s>0$, and that
$\max(\psi(-t),\psi(t))\le \widehat{\psi}(t)\le\psi(t)+\psi(-t)$ for all
$t\in\Bbb{R}$.  From the fixed point equation $\zi\od Q+M\zi$, we infer
$\psi(t)=\me e^{tQ}\psi(Mt)$ for all $t\in\Bbb{R}$. These facts will be used in several
places hereafter.
\medskip

(a) Sufficiency: Suppose that (\ref{cond4000}) is valid.
The almost sure finiteness of $\zi$ follows from Proposition \ref{exper}. We have to
check that $r(\zi)\geq r(Q)$. To this end, we fix an arbitrary $s\in (0,r(Q))$ and
divide the subsequent proof into two steps.
\smallskip

Step (a1). Assume first that $|M|\leq \beta<1$ a.s. for some $\beta>0$. Since the
function $\widehat\varphi$ is convex and differentiable on
$[0,\beta s]$, its derivative is nondecreasing on that interval. Therefore,
for each $k=2,3,\ldots$, there exists $\theta_{k}\in
[0,\beta^{k-1}s]$ such that
\begin{equation}\label{alls}
0\ \le\ \widehat{\varphi}(\beta^{k-1}s)-1\ =\ \widehat{\varphi}'(\theta_{k})
\beta^{k-1}s\ \le\ \widehat{\varphi}^\prime (\beta s)\beta^{k-1}s.
\end{equation}
Now $r(\zi)\geq r(Q)$ follows from
\begin{eqnarray*}
\widehat{\psi}(s) &\leq& \me\exp\left(s\sum_{k\ge 1}|\Pi_{k-1}Q_k|\right)\\
&\le&\me\exp\left(s\sum_{k\ge 1}\beta^{k-1}|Q_k|\right)
\ =\ \prod_{k\ge 1}\widehat{\varphi}(\beta^{k-1}s)\\
&\leq& \widehat{\varphi}(s)\exp\left(\sum_{k\ge 2}
(\widehat{\varphi}(\beta^{k-1}s)-1)\right)\\
&\le&\widehat{\varphi}(s)\exp\Big(\widehat{\varphi}^\prime
(\beta s)\beta s (1-\beta)^{-1}\Big)\ <\ \infty.\qquad\text{[by (\ref{alls})]}
\end{eqnarray*}

Step (a2). Consider now the general case. Since $\mmp\{|M|=1\}=0$,
we can choose $\beta\in (0,1)$ such that
$$\mmp\{|M|>\beta\}\ <\ 1\quad\text{and}\quad\gamma\ \gl\ \me e^{s|Q|}
\1_{\{|M|>\beta\}}<1. $$
Define the a.s.\ finite stopping times
$$\ T_0\ \gl\ 0,\quad T_k\ \gl\ \inf\{n>T_{k-1}: |M_n|\leq \beta\},\quad k=1,2,\ldots $$
We have $\zi=Q_1^\ast+\sum_{k\ge 1}
M_1^\ast\cdot...\cdot M_{k-1}^\ast Q_k^\ast$, where for $k=1,2,\ldots$
\begin{eqnarray}
M_k^\ast &\gl& M_{T_{k-1}+1} \cdot...\cdot M_{T_k}\ =\ \Pi_{T_{k-1}+1:T_{k}}
\quad\text{and}\label{def1}\\
Q_k^\ast &\gl& Q_{T_{k-1}+1}+M_{T_{k-1}+1}Q_{T_{k-1}+2}+\ldots+M_{T_{k-1}+1}\cdot...
\cdot M_{T_k-1}Q_{T_k},\qquad\label{def2}
\end{eqnarray}
so that $(M_k^\ast, Q_k^\ast)$ are independent copies of
$$ (M^\ast, Q^\ast)\ \gl\ \big(\Pi_{T_1},
Q_1+\textstyle\sum_{k=1}^{T_1}\Pi_{k-1}Q_k\big). $$
Since $|M^\ast|\leq\beta$ a.s., Step (a1) of the proof provides the desired conclusion
if we still verify that $\widehat{\varphi}(s)<\infty$ implies
$\me e^{s|Q^\ast|}<\infty$. This is checked as follows:
\begin{eqnarray*}
\me e^{s|Q^\ast|} &\leq& \me e^{s(|Q_1|+\ldots+|Q_{T_1}|)}
\ =\ \sum_{n\ge 1}\me e^{s(|Q_1|+\ldots+|Q_n|)}\1_{\{T_1=n\}}\\
&=&\sum_{n\ge 1} \me\Bigg[e^{s|Q_n|}\1_{\{|M_n|\leq\beta\}}\prod_{k=1}^{n-1}
e^{s|Q_{k}|}\1_{\{|M_{k}|>\beta\}}\Bigg]\\
&=&\me e^{s|Q|}\1_{\{|M|\leq \beta\}}\sum_{n\ge 1}\gamma^{n-1}\\
&\leq&\widehat{\varphi}(s)(1-\gamma)^{-1}\ <\ \infty.
\end{eqnarray*}

(b) Necessity: If $\me e^{s|\zi|}<\infty$, we have $\me |\zi|^p<\infty$ and therefore,
by Theorem \ref{mom1},
$\me|M|^{p}<1$ for all $p>0$. The latter in combination with $\mmp\{|M|=1\}=0$ implies
$|M|<1$ a.s. Finally, if
$\widehat{\psi} (s)<\infty$ and $c\gl\min_{|t|\le s}\psi(t)$ (clearly $>0$), then
\begin{equation}\label{qmoment}
\infty\ >\ \psi(t)\ =\ \me e^{tQ}\psi(Mt)\ \ge\ c\varphi(t),\quad t\in\{-s,s\},
\end{equation}
and thus $\widehat{\varphi}(s)\le\varphi(s)+\varphi(-s)<\infty$. This shows $r(\zi)\le
r(Q)$.
\end{myproof}
\bigskip

\begin{myproof}{\emph{of Theorem \ref{perexp2}}}. Recall that
$a_{\pm}\gl\me e^{\pm sQ}\1_{\{M=1\}}$, $b_{\pm}\gl\me e^{\pm sQ}\1_{\{M=-1\}}$.

(a) Necessity:
By the same argument as in part (b) of the proof of Theorem \ref{perexp}, we infer
$\me|M|^{p}<1$ for all $p>0$ and thereby $|M|\le 1$ a.s. Moreover, as $\widehat{\psi}
(s)<\infty$, inequality (\ref{qmoment}) holds here as well and gives $\widehat
{\varphi}(s)<\infty$. This shows (\ref{cond4001}) and leaves us with the proof of
(\ref{cond4002}), for
which we will proceed in two steps:
\smallskip

Step (a1). Suppose first that $\mmp\{M=-1\}=0$ in which case $b_{\pm}=0$ and
thus (\ref{cond4002}) reduces to $a_{\pm}<1$. We have
\begin{eqnarray*}
\psi(s)&=&\me e^{sQ}\psi(Ms)\1_{\{|M|<1\}}+\psi(s)\,\me e^{sQ}
\1_{\{M=1\}}\quad\text{and}\\
\psi(-s)&=&\me e^{-sQ}\psi(-Ms)\1_{\{|M|<1\}}+\psi(-s)\,\me e^{-sQ}\1_{\{M=1\}},
\end{eqnarray*}
which together with $\me e^{sQ}\psi(\pm Ms)\1_{\{|M|<1\}}>0$ (as $\mmp\{|M|<1\}>0$)
implies
$$ \me e^{\pm sQ}\1_{\{M=1\}}=a_{\pm}<1 $$
as required.
\medskip

Step (a2). Assuming now $\mmp\{M=-1\}>0$, let
$\{(M_k^{*},Q_k^{*}):k=1,2,...\}$ be defined as in (\ref{def1}),\,(\ref{def2}), but with
$$ T_0\ \gl\ 0,\quad T_k\ \gl\ \inf\{n>T_{k-1}: \Pi_{T_{k-1}+1:n}>-1\},
\quad k=1,2,\ldots $$
Then $\mmp\{M^{*}=-1\}=0$, and we infer from Step (a1) that $\me e^{\pm sQ^{*}}
\1_{\{M^{*}=1\}}<1$. But
\begin{eqnarray*}
e^{\pm sQ_{1}^{*}}\1_{\{M_{1}^{*}=1\}}&=&e^{\pm sQ_{1}}\1_{\{M_{1}=1\}}\ +
\ e^{\pm s(Q_{1}-Q_{2})}\1_{\{M_{1}=M_{2}=-1\}}\\
&+&\sum_{n\ge 3}e^{\pm s(Q_{1}-Q_{2}-...-Q_{n})}\1_{\{M_{1}=-1,M_{2}=...=M_{n-1}=1,
M_{n}=-1\}}
\end{eqnarray*}
implies
\begin{equation*}\label{equone}
1\ >\ \me e^{\pm sQ^{*}}\1_{\{M^{*}=1\}}\ =\ a_{\pm}+\sum_{n\ge 0}b_{\pm}
a_{\mp}^{n}b_{\mp}\ =\ a_{\pm}+\dfrac{b_{\pm}b_{\mp}}{1-a_{\mp}}
\end{equation*}
and thus (\ref{cond4002}).

(b) Sufficiency: Let $\{(M_k^\ast, Q_k^\ast):\,k=1,2,...\}$ be as defined in Step
(a2). Assuming (\ref{cond4002}) we thus have
$$ a_{\pm}^{*}\ \gl\ \me e^{\pm sQ^{*}}\1_{\{M^{*}=1\}}\ =\ a_{\pm}+\dfrac{b_{\pm}
b_{\mp}}{1-a_{\mp}}\ <\ 1. $$
Note that (\ref{cond4002}) particularly implies $a_{\pm},b_{\pm}\in [0,1)$. Using this
and
\begin{eqnarray*}
e^{\pm sQ_{1}^{*}}&=&e^{\pm sQ_{1}}\1_{\{M_{1}>-1\}}\ +\ e^{\pm s(Q_{1}-Q_{2})}
\1_{\{M_{1}=-1,M_{2}>-1\}}\\
&+&\sum_{n\ge 3}e^{\pm s(Q_{1}-Q_{2}-...-Q_{n})}\1_{\{M_{1}=-1,M_{2}=...=M_{n-1}=1,
M_{n}>-1\}}
\end{eqnarray*}
we further obtain that
$$ \me e^{\pm sQ^{*}}\ =\ \me e^{\pm sQ}\1_{\{M>-1\}}\ +\ \me e^{\mp sQ}\1_{\{M>-1\}}
{b_{\pm}\over 1-a_{\mp}}\ <\ \infty. $$
Now let
$$ \widehat{T}_0\ \gl\ 0,\quad\widehat{T}_k\ \gl\ \inf\{n>T_{k-1}:
M_{n}^{*}<1\},\quad k=1,2,\ldots $$
and then $\{(\widehat{M}_k,\widehat{Q}_k):k=1,2,...\}$ in accordance with
(\ref{def1}),\,(\ref{def2}) for these stopping times. We claim that $\me e^{\pm s
\widehat{Q}}<\infty$ and thus $\me e^{s|\widehat{Q}|}<\infty$. Indeed,
\begin{eqnarray*}
\me e^{\pm s\widehat{Q}}&=&\sum_{n\ge 1}\me e^{\pm s(Q_{1}^{*}+M_{1}^{*}Q_{2}^{*}
+...+M_{1}^{*}\cdot...\cdot M_{n-1}^{*}Q_{n}^{*})}
\1_{\{M_{1}^{*}=...=M_{n-1}^{*}=1,M_{n}^{*}<1\}}\\
&=&\me e^{\pm sQ^{*}}\1_{\{M^{*}<1\}}\sum_{n\ge 1}a_{\pm}^{*n-1}\\
&\le&\dfrac{\me e^{\pm sQ^{*}}}{1-a_{\pm}}\ <\ \infty.
\end{eqnarray*}
So we have $\mmp\{|\widehat{M}|=1\}=0$ and $\me e^{s|\widehat{Q}|}<\infty$ and may thus
invoke Theorem \ref{perexp} to finally conclude $\me e^{s|\zi|}<\infty$ because
$\zi$ is also the perpetuity generated by $(\widehat{M},\widehat{Q})$. \end{myproof}

\section{Examples}

We begin with an example that shows that condition
$\mmp\{M=0\}=0$ in Theorem \ref{continuity} is indispensable.
\begin{example}\rm
If $\mmp\{M=0\}=p=1-\mmp\{M=1\}$ for some $p\in (0,1)$ and $Q=1$ a.s., the
distribution of $\zi$ is geometric with parameter $p$. This can be seen
from Remark \ref{geom}, as the random variable $N$ defined there
has a geometric distribution with parameter $\mmp\{M=0\}$.
\end{example}
The next examples illustrate that the distribution of $\zi$ can indeed be
continuously singular as well as absolutely continuous. Denote by $\mL(X)$ the
distribution of a random variable $X$.
\begin{example}
\label{singSN} \rm [Deterministic $M$] Consider the situation where $M$ is a.s.\
equal to a constant $c\in (0,1)$, so
$$\zi \od c\zi+Q. $$

(a) If $c=1/2$ and $\mL(Q)$ is a Poisson distribution with parameter $\lambda>0$,
then $\mL(\zi)$ is singularly continuous according to Example 4.3 by Watanabe
\cite{Watanabe}.

(b) If $c=1/n$ for some fixed positive integer $n$ and $\mL(Q)$ is the discrete
uniform distribution on $\{0,...,n-1\}$, then it can be easily verified that
$\zi$ has the uniform distribution on $(0,1)$ and is thus absolutely continuous.
This example is a special case of one due to Letac, see Example A2 in \cite{Let}.

(c) A particularly well-studied class of special cases is when
$\mmp\{Q=1\}=\mmp\{Q=-1\}=1/2$. A short survey can be found in \cite{Diaconis}.
If $c=1/2$, then $\zi$ is uniformly distributed on $[-2,2]$, while $\zi$ is
continuously singular if $0<c<1/2$. One would expect $\zi$ to be absolutely continuous
whenever $c\in (1/2,1)$. However, this is not true as there are values
of $c$ between $1/2$ and 1 giving a singular $\mL(\zi)$, for example, if
$c=(\sqrt{5}-1)/2=0.618...$, see \cite{Erdos1}, \cite{Erdos2}. Meanwhile it
has been proved by Solomyak \cite{Solomyak} that, on the other hand, $\mL(\zi)$
is indeed absolutely continuous for almost all values of $c\in (1/2,1)$.
\end{example}

\begin{example}\label{PYgamma}\rm Assume that $M$ and $Q$ are independent.
\medskip

(a) If $M$ has a beta distribution with parameters 1 and $\alpha>0$, i.e.
$$ \mmp\{M\in dx\}\ =\ \alpha(1-x)^{\alpha-1}\1_{(0,1)}(x)\ dx, $$
and if $Q$ has a $\Gamma(\alpha,\alpha)$-distribution, i.e.
$$ \mmp\{Q\in dx\}\ =\ \dfrac{\alpha^\alpha}{\Gamma(\alpha)}x^{\alpha-1}e^{-\alpha
x}\1_{(0,\infty)}(x)\ dx, $$
then $\mL(\zi)=\Gamma(\alpha+1,\alpha)$ as one can easily verify by direct calculation.
\medskip

(b) If $M$ has a Weibull distribution with parameter $1/2$, i.e.
$$ \mmp\{M\in dx\}\ =\ \dfrac{e^{-\sqrt{x}}}{2\sqrt{x}}\1_{(0,\infty)}(x)\ dx, $$
and $Q$ is nonnegative with Laplace transform $\me e^{-sQ}=(1+b\sqrt{s})
e^{-b\sqrt{s}}$,
$s\geq 0$, $b>0$, then
$$\int_0^\infty e^{-sx}\mmp\{\zi\in dx\}\ =\ e^{-b\sqrt{s}},\quad s\geq 0$$
i.e. $\mL(\zi)$ is the positive stable law with index $1/2$, as was found independently
in \cite{Donati} and \cite{IksKim2}.
\medskip

(c) If $\mmp\{M\in dx\}=(x^{-1/2}-1)\1_{(0,1)}(x)dx$ and $Q$ has Laplace transform
$$ \varphi(s)\ \gl\ \left(\dfrac{\sqrt{2s}}{\sinh \sqrt{2s}}\right)^2,\quad
s\geq 0, $$
then the Laplace transform of $\zi$ takes the form
$$ \me e^{-s\zi}\ =\ \dfrac{3(\sinh \sqrt{2s}-\sqrt{2s}\cosh \sqrt{2s})}{\sinh^3
\sqrt{2s}}\ =\ -\dfrac{\varphi'(s)}{\me Q},\quad s\geq 0. $$
This result was obtained in \cite{PitmanYor}.
\end{example}
\bigskip

Absolute continuity of $\mL(\zi)$ in (a) and (b) is obvious. In (c) it
follows from the fact that the corresponding $\mL(M)$ is absolutely
continuous.

Recall that the \emph{size-biased distribution} $\omu$ pertaining to a probability
distribution $\mu$ on $[0,\infty)$ with finite mean $m>0$ is defined as
$$ \omu(dx)\ \gl\ m^{-1}x\ \mu (dx). $$
In all three previous examples $\mL(\zi)$ is the size-biased distribution pertaining
to $\mL(Q)$. The study of distributions solving the fixed point equation
(\ref{fpeq}) and having this additional property was initiated by Pitman and Yor
\cite{PitmanYor} and then continued in \cite{IksKim1} and \cite{IksKim2}.

Our last example provides an illustration of Theorem \ref{perexp2}.

\begin{example}\rm
Let $Q$ be an exponential random variable with parameter $a>0$ and
$M$ be independent of $Q$ with $\mmp\{0\le M\le 1\}=1$ and $\me M<1$.
Then $\zi$ is a.s.\ finite and it can be checked directly or by
using the fact that $\zi \od \int_0^\infty e^{-Y_t}dt$ for an
appropriate compound Poisson process $\{Y_{t}:t\ge 0\}$ starting at zero, that
$$ \me \zi^n\ =\ \dfrac{n!}{a^n (1-\me M)(1-\me M^2)\cdots (1-\me
M^n)}. $$
Put $a_n\gl\me\zi^n/n!$ and note that $\lin a_{n+1}^{-1}a_{n}=a\,\mmp\{M<1\}$.
Hence, by the Cauchy-Hadamard formula,
$r(\zi)=a \mmp\{M<1\}$ which is in full accordance with Theorem
\ref{perexp2} according to which $r(\zi)$ is the positive solution to the equation
$\frac{a}{a-s}\mmp\{M=1\}=1$ and thus indeed equal to $a\mmp\{M<1\}$.
\end{example}

\section{Concluding remarks}

Although settling a number of open questions about perpetuities, this work gives also
rise for further research. For example, it is natural to ask for conditions under which
$\zi$ is of each of the three possible types. As already pointed out after Theorem \ref{continuity}, the law of $\zi$ can only be degenerate if $\mmp\{Q+cM=c\}=1$
for some $c\in\Bbb{R}$. However, the discussion in Example \ref{singSN}(c) indicates that a similar characterization for singularity and absolute continuity of the law of $\zi$ remains open even in the special case where $M$ is deterministic and $Q$ very simple.

Another natural problem that arises when regarding our Theorems
\ref{perexp} and \ref{perexp2} is to determine $r_{*}(\zi)\le 0$ and $r^{*}(\zi)\ge 0$,
given by
$$ r_{*}(\zi)\gl\inf\{r\le 0:\me e^{s\zi}<\infty\text{ for all }r\le s\le 0\} $$
and
$$ r^{*}(\zi)\gl\sup\{r\ge 0:\me e^{s\zi}<\infty\text{ for all }0\le s\le r\}. $$
While this would provide information about when $\zi$ has exponential left and right
tails, one may also aim at conditions that ensure existence of log-type moments
of $\zi$ of the form $\me(\log^{+}|\zi|)^{\beta}$ for $\beta>0$. The latter is studied in
recent work by the first two authors \cite{AlsIks} and is of additional interest due to
the connection of perpetuities with certain intrinsic martingales in the
supercritical branching random walk. This connection was first observed in \cite{Iks041}
and further exploited in \cite{IksRos}.

\bigskip

\textbf{Acknowledgment.} The main part of this work was done while
A.\ Iksanov was visiting the Institute of Mathematical Statistics
at M\"{u}nster in October/November 2006. He gratefully
acknowledges financial support and hospitality. The authors are
further indebted to an Associate Editor and two anonymous referees
for making a number of valuable remarks that helped improving the
presentation of this work. The research of A.\ Iksanov and U.\
R\"{o}sler was supported by the DFG grant, project no.436UKR
113/93/0-1.

\end{document}